\documentclass[10pt]{article}
\usepackage{epsfig}
\usepackage{graphicx,graphics}
\usepackage{pictex}
\usepackage{amsmath,amssymb}
\usepackage{textcomp}
\usepackage{verbatim}

\newtheorem{thm}{Theorem}[section]
\newtheorem{prop}[thm]{Proposition}
\newtheorem{cor}[thm]{Corollary}
\newtheorem{lem}[thm]{Lemma}
\newtheorem{conj}[thm]{Conjecture}
\newtheorem{que}[thm]{Question}
\newtheorem{exa}[thm]{Example}
\newtheorem{defn}[thm]{Definition}
\newcommand{\ben}{\begin{enumerate}}
\newcommand{\een}{\end{enumerate}}
\newcommand{\ble}{\begin{lem}}
\newcommand{\ele}{\end{lem}}
\newcommand{\bthm}{\begin{thm}}
\newcommand{\ethm}{\end{thm}}
\newcommand{\bpr}{\begin{prop}}
\newcommand{\epr}{\end{prop}}
\newcommand{\bco}{\begin{cor}}
\newcommand{\eco}{\end{cor}}
\newcommand{\bcon}{\begin{conj}}
\newcommand{\econ}{\end{conj}}
\newcommand{\bque}{\begin{que}}
\newcommand{\eque}{\end{que}}
\newcommand{\bde}{\begin{defn}}
\newcommand{\ede}{\end{defn}}
\newcommand{\bex}{\begin{exa}}
\newcommand{\eex}{\end{exa}}
\newcommand{\barr}{\begin{array}}
\newcommand{\earr}{\end{array}}
\newcommand{\btab}{\begin{tabular}}
\newcommand{\etab}{\end{tabular}}
\newcommand{\beq}{\begin{equation}}
\newcommand{\eeq}{\end{equation}}
\newcommand{\bea}{\begin{eqnarray*}}
\newcommand{\eea}{\end{eqnarray*}}
\newcommand{\bce}{\begin{center}}
\newcommand{\ece}{\end{center}}

\newtheorem{obs}[thm]{Observation}
\newcommand{\bobs}{\begin{obs}}
\newcommand{\eobs}{\end{obs}}
\newtheorem{prob}[thm]{Problem}
\newcommand{\bprob}{\begin{prob}}
\newcommand{\eprob}{\end{prob}}
\newcommand{\pf}{{\bf Proof.}}

\newcommand{\qed}{\rule{1ex}{1ex}}

\newcommand{\ol}{\overline}

\newcommand{\sbe}{\subseteq}

\newcommand{\de}{\delta}

\newcommand{\cP}{{\cal P}}

\newcommand{\diam}{\mathop{\rm diam}}

\newcommand{\corona}{\mathop{\rm cor}}

\newcommand{\maj}{\mathop{\rm maj}}
\newcommand{\smaj}{\mathop{\rm smaj}}

\newcommand{\roy}{\mathop{\rm roy}}
\newcommand{\sroy}{\mathop{\rm sroy}}

\usepackage{amsmath}

\begin{document}

\title{Extremal Problems in\\
Royal Colorings of Graphs}

\author{$^{1}$Akbar Ali,  $^{2}$Gary Chartrand,\\
$^{2}$James Hallas and $^2$Ping Zhang\\
\\
$^1$University of Management and Technology\\
                Sialkot 51310, Pakistan\\
$^1$College of Science, University of Hail\\
Hail 81451, Saudi Arabia\\
\\
$^2$Western Michigan University\\
Kalamazoo, Michigan 49008, USA\\
Email: ping.zhang@wmich.edu}

\date{}

\maketitle

\bce
{\bf Abstract}
\ece

\begin{quote}
An edge coloring~$c$ of a graph~$G$ is a royal $k$-edge coloring of~$G$ if the edges of~$G$ are assigned
nonempty subsets of the set $\{1, 2, \ldots, k\}$ in such a way that the vertex coloring
obtained by assigning the union of the colors of the incident edges of each vertex
is a proper vertex coloring. If the vertex coloring is  vertex-distinguishing, then $c$ is   a   strong  royal  $k$-edge coloring.
 The minimum positive integer~$k$ for which~$G$   has a strong royal $k$-edge coloring is
the   strong royal index  of~$G$. It has been conjectured that if
$G$  is a  connected graph of order $n\ge 4$ where $2^{k-1} \le n \le 2^k-1$ for a positive integer~$k$, then
the  strong royal index  of~$G$ is either~$k$ or~$k+1$.  We discuss this conjecture along with other information
concerning strong royal colorings of graphs.  A sufficient condition for such a graph to have strong royal index~$k+1$ is presented.
  \end{quote}

\noindent
{\bf Key Words:}   color-induced  coloring,  royal and strong royal coloring,  strong royal  index.

%\vs{.1cm}

\noindent
{\bf AMS Subject Classification:} 05C15, 05C05, 05C35.

%\doublespacing

\section{Introduction}

During the past several years, a number of edge colorings (or  edge labelings) have been introduced that
 give rise to  vertex colorings  that are either proper  or  vertex-distinguishing (see~\cite{att, Burris95, Burris97, cjlors88}, for example).
 Many of these are discussed in the books~\cite{HLG, CGT}.  We discuss two of these colorings here.
 For a connected graph $G$ of order~3 or more and a positive integer~$k$,
let  $c: E(G) \to [k]=\{1, 2, \ldots, k\}$ be an unrestricted edge coloring of~$G$, that is, adjacent edges of~$G$ may be assigned the same color.
We write $\cP^*([k])$ for the set consisting of the $2^k-1$ nonempty subsets of~$[k]$.
The edge coloring $c$ gives rise to the   vertex coloring $c': V(G) \to \cP^*([k])$
where  $c'(v)$ is the set of  colors of the edges incident with $v$.
If  $c'$  is a proper vertex coloring of~$G$, then  $c$ is a  {\it majestic   $k$-edge coloring}
and the minimum positive integer~$k$ for which $G$ has a majestic   $k$-edge coloring
is the  {\it  majestic   index} $\maj(G)$ of~$G$.
If $c'$   is  {\it vertex-distinguishing} (that is, $c'(u)\ne c'(v)$ for every two distinct
vertices $u$ and $v$ of~$G$),  then  $c$ is a  {\it strong majestic   $k$-edge coloring}
and the minimum positive integer~$k$ for which $G$ has a strong majestic   $k$-edge coloring
is the  {\it  strong majestic   index} $\smaj(G)$ of~$G$.
Majestic  edge colorings  were  introduced by   Gy\"{o}ri,     Hor\v{n}\'{a}k,    Palmer, and   Wo\'{z}nick~\cite{gnid1}
under different terminology and studied further  in~\cite{IHart2018, gnid2}.
Strong  majestic edge colorings  were introduced by  Harary and  Plantholt~\cite{Harary85} in 1985, also using different terminology,  and
studied further by others (see~\cite{CGT, zhang, zhang2}).

 While  an edge coloring $c$ of a graph~$G$ typically  uses  colors from  the set $[k]$ for some positive integer~$k$ resulting in
$c(e)=i$  for some $i \in [k]$, we might equivalently  define $c(e)=\{i\}$ as well. Expressing the edge coloring~$c$ in this way results in
both~$c$ and the induced vertex coloring $c'$ assigning subsets of~$[k]$ to   the edges as well as the vertices  of~$G$.
Furthermore, expressing~$c$ in this manner  suggests the idea of studying edge colorings $c$ where both $c$ and its derived vertex coloring~$c'$
assign nonempty subsets of~$[k]$ to the elements (edges and vertices) of a graph~$G$ such that the color assigned to an edge of~$G$ by~$c$
is not necessarily a singleton subset of~$[k]$.   This observation gives rise to the primary concepts of this paper,
namely royal and strong royal colorings,  which were introduced in~\cite{royc1}.

For a   positive integer~$k$,   let $\cP^*([k])$ denote the collection of the $2^k-1$ nonempty subsets of the set~$[k]$.
For a connected graph~$G$ of order~3 or more, an edge coloring $c: E(G) \to \cP^*([k])$ of~$G$
is a {\it royal $k$-edge coloring}  if  the vertex  coloring $c': V(G) \to \cP^*([k])$  defined by
                  $c'(v)=\bigcup_{e \in E_v} c(e),$
where $E_v$ is the set of edges of~$G$ incident with~$v$,  is   proper, that is,
adjacent vertices are assigned distinct colors.
If $c'$ is vertex-distinguishing, then  $c$ is   a  {\it  strong  royal $k$-edge coloring} of~$G$.
The minimum positive integer~$k$ for which~$G$ has a strong royal $k$-edge coloring is
the {\it  strong royal index} of~$G$, denoted by~$\sroy(G)$.   Therefore,
royal colorings  are generalizations of majestic  edge colorings and
strong royal colorings  are generalizations of  strong majestic   colorings.
This concept  was independently introduced and studied  in~\cite{BDK2017, royc1}.
While there are many connected graphs~$G$ for which $\sroy(G) \ne  \smaj(G)$,
we know of no graph~$G$ for which $\roy(G) \ne  \maj(G)$.
Consequently, our emphasis  here is on the strong royal indexes of graphs.
If $G$  is a  connected graph of order $n \ge 4$,  there is a unique integer $k \ge 3$  such that   $2^{k-1} \le n \le 2^k-1$.
We now present several useful  observations  made in~\cite{BDK2017, royc1}.

\bobs        \label{sroyuppern}
\   If $G$ is a connected graph of order~$n \ge 4$  where $2^{k-1} \le n \le 2^k-1$, then
%  \bce
    $\sroy(G) \ge k.$
%  \ece
\eobs

\bobs            \label{royalspanH}
\  If  $G$ is a connected graph of order~$4$ or more, then
    $$\sroy(G) \le 1+ \min\{\sroy(H):  \mbox{ $H$ is a connected spanning subgraph of~$G$}\}.$$
In particular,  $\sroy(G) \le 1+ \min\{\sroy(T):  \mbox{ $T$ is a spanning tree of~$G$}\}.$
 \eobs

It was shown in~\cite{BDK2017} that if $G$ is a connected graph of order~$n \ge 4$  where $2^{k-1} \le n \le 2^k-1$, then $\sroy(G) \le k+2$.
Furthermore, it was conjectured  in~\cite{royc1} that the strong royal index of every connected graph of   order $n \ge 4$ where $2^{k-1} \le n \le 2^k-1$ is
either~$k$ or~$k+1$. This gives rise to the following concepts.
%
%\bcon     \label{royalgraph2}
%\   If $G$  is a  connected graph of order $n \ge 4$ where $2^{k-1} \le n \le 2^k-1$ for some integer~$k$, then
%$\sroy(G)$ is either~$k$ or~$k+1$.
%\econ
%
A   connected graph~$G$  of order $n \ge 3$ where $2^{k-1} \le n \le 2^k-1$ is a {\it  royal-zero graph}  if $\sroy(G)=k$
and is a {\it  royal-one  graph} if $\sroy(G)=k+1$.  Therefore, the conjecture on the strong royal index can  be rephrased as follows.
%Conjecture~\ref{royalgraph2} can then be rephrased as follows.

 \bcon                      \label{royalgraph2.1}
\  Every connected graph  of order at least~$4$ is either  royal-zero or  royal-one.
\econ

By Observation~\ref{royalspanH},   the strong royal indexes of  trees play an important role  in the study of  strong royal indexes of   connected graphs.
It was conjectured in~\cite{royc1} that every tree  of  order $n \ge 4$ where $2^{k-1} \le n \le 2^k-1$ has strong royal index~$k$ and consequently is royal-zero.
This conjecture can therefore be  rephrased in terms of royal-zero graphs.

\bcon                   \label{royaltree}
\   Every tree  of order at least~$4$ is royal-zero.
\econ

%\bcon                   \label{royaltree}
%\   If $T$  is a tree of order $n \ge 4$ where $2^{k-1} \le n \le 2^k-1$ for some integer~$k$, then
%$\sroy(T) = k$ and so $T$ is royal-zero.
%\econ

Conjecture~\ref{royaltree}  has been verified for trees of small order (order~10 or less), all paths,
all complete binary trees, all caterpillars of diameter~4 or less
as well as some specialized trees (see~\cite{BDK2017, royc1}).
%(A {\it caterpillar} is a tree~$T$ of order 3 or more, the removal of whose leaves  produces a path.)
By Observation~\ref{royalspanH}, it follows that if Conjecture~\ref{royaltree} is true, then
Conjecture~\ref{royalgraph2.1} is true as well.
 While the strong royal index of each cycle was stated in~\cite{BDK2017},
 we illustrate the  concepts described above by  providing a proof that describes in each case an appropriate
 edge coloring.

\begin{figure}[ht]\centering
\centerline{\input{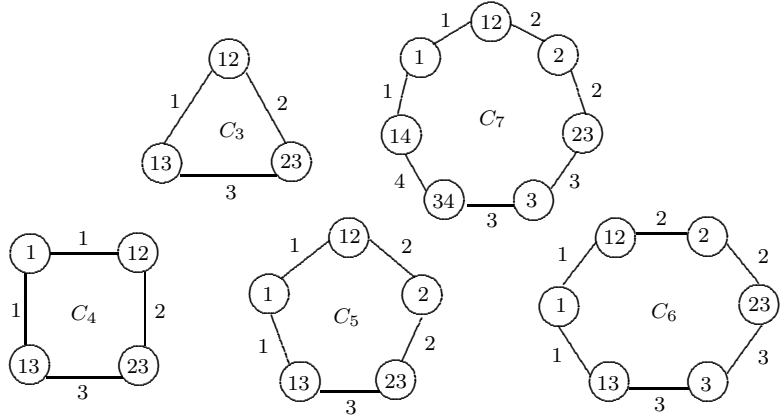}}
\caption{Strong royal colorings of $C_n$  where $3 \le n \le 7$} \label{roy21}
\end{figure}

\bex                \label{royalcn}
\   For every  integer $n \ge 3$,
     $$\sroy(C_n) =\left\{\barr{cl}
   1+  \left\lceil{\log_2 (n+1)}\right\rceil   & \mbox{ if $n = 3, 7$}\\[.2cm]
     \left\lceil{\log_2 (n+1)}\right\rceil  & \mbox{ if $n \ne  3, 7$.}
     \earr\right.
     $$
That is,  if  $C_n$  is a  cycle  of length $n \ge 3$ where $2^{k-1} \le n \le 2^k-1$ for some integer~$k$, then
$\sroy(C_n) = k$  unless $n = 3$ or $n = 7$, in which case, $\sroy(C_3)= 3$ and $\sroy(C_7)=4$.
 \eex
\pf \quad  Let $k =  \left\lceil{\log_2 (n+1)}\right\rceil \ge 2$. Then $2^{k-1} \le n \le 2^k-1$.
We show that $\sroy(C_3) = 3$, $\sroy(C_7) = 4$, and $\sroy(C_n) = k$ if $n \ne 3, 7$.
Figure~\ref{roy21} shows a strong royal 3-edge coloring  of $C_3$ and a strong royal 4-edge coloring  of~$C_7$,
which shows that $\sroy(C_3) \le 3$ and $\sroy(C_7) \le 4$.
(For simplicity, we write  the set $\{a\}$ as $a$, $\{a, b\}$ as $ab$,   and $\{a, b, c\}$ as $abc$.)
If $\sroy(C_3)=2$, then because $|\cP^*([2])| = 3$, there are vertices of~$C_3$ colored~1 and~2, implying  that two edges
of~$C_3$ are colored with each of these two colors, which is impossible.
If $\sroy(C_7)=3$, then because $|\cP^*([3])| = 7$, there are vertices of~$C_7$ colored~1, 2,  and~3, implying that two edges
of~$C_7$ are colored with each of these three colors.   Regardless of how the seventh edge of~$C_7$ is colored, the resulting set of vertex colors
is not~$\cP^*([3])$. Consequently,  $\sroy(C_3) = 3$ and $\sroy(C_7)= 4$.  By Proposition~\ref {sroyuppern}, it suffices to show
that $C_n$ has a strong  royal $k$-edge coloring if $n \ne 3, 7$.
Figure~\ref{roy21}  also shows a strong royal 3-edge coloring for each of $C_4, C_5$, and $C_6$
and so $\sroy(C_n) = 3$  for $n = 4, 5, 6$.

Next, suppose  that $n \ge 8$, where $2^{k-1} \le n \le 2^k-1$
for a unique integer $k \ge 4$. We show that  $C_n$ has a strong  royal $k$-edge coloring
by  considering two cases, depending on whether $n$ is even or $n$ is odd.
Let $P_n=(v_1,v_2,\ldots,v_n)$ where  $e_i=v_iv_{i+1}$ for $1\leq i\leq n-1$.

\smallskip

\begin{figure}[ht]\centering
\centerline{\input{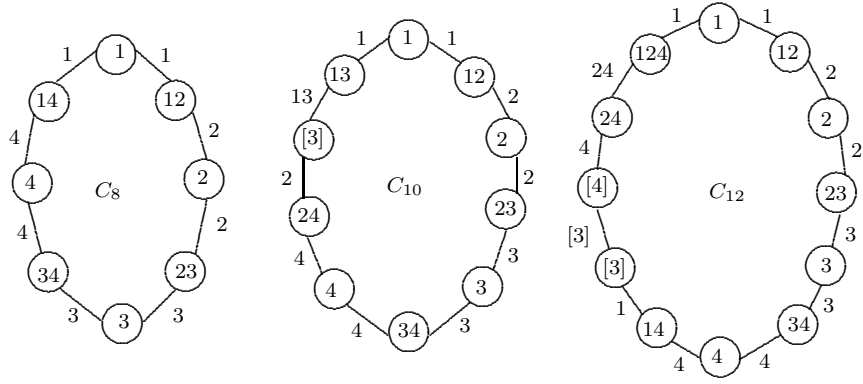}}
\caption{Strong royal $4$-edge colorings of $C_n$  for $n =8, 10, 12$} \label{roy24}
\end{figure}

{\it Case~$1$.} {\it $n \ge 8$ is even}. Figure~\ref{roy24}   shows a strong royal 4-edge coloring for each of~$C_8$, $C_{10}$, and~$C_{12}$
and so $\sroy(C_n) = 4$  for $n =8, 10, 12$.

\begin{figure}[ht]\centering
\includegraphics[width=10cm]{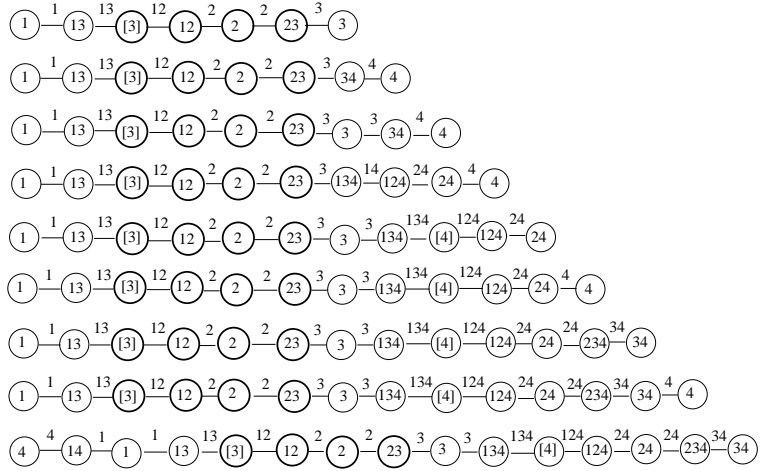}
%\centerline{\input{roy27.pictex}}
\caption{Strong royal $(k-1)$-edge colorings of $P_r$  for $7\le r \le 15$} \label{roy27}
\end{figure}

Thus, we assume that  $n = 2r \ge 14$ where $r \ge 7$ is an integer such that $2^{k-2} \le r \le 2^{k-1}-1$.
If  $r = 7$, then $k-1=3$; while if $8 \le r \le 15$, then $k-1=4$.
A  strong royal $(k-1)$-edge coloring  $c$ for each path $P_r$ ($7 \le r \le 15$) is shown  in Figure~\ref{roy27}.

\begin{figure}[ht]\centering
\centerline{\input{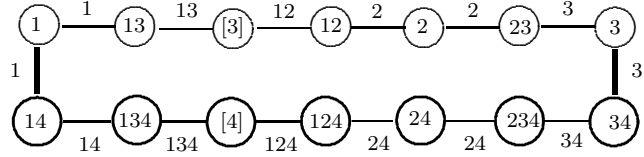}}
\caption{Constructing a strong royal $4$-edge coloring  of $C_{14}$} \label{roy25}
\end{figure}

For $7 \le r \le 15$, let $P_r=(v_1, v_2, \ldots, v_r)$ and let $P^*_r=(u_1, u_1, \ldots, u_r)$.
The path $P_{2r}$ is constructed from $P_r$ and $P^*_r$ by adding the edge $v_ru_r$ and the cycle
$C_{2r}$ is constructed from $P_{2r}$ by adding the edge $v_1u_1$. The  edge coloring $c$ is extended
(1) to an edge coloring~$c$ of~$P_{2r}$ by defining $c(u_iu_{i+1}) = c(v_iv_{i+1}) \cup \{k\}$
        (where $k = 4$ if $r = 7$ and $k = 5$ if $8 \le r \le 15$) for $1 \le i \le r-1$ and  $c(v_ru_r)= c(v_{r-1}v_r)$ and
(2) to an edge coloring~$c$ of~$C_{2r}$ by defining  $c(v_1u_1)=c(v_1v_2)$ in addition.
In this manner, no vertex of~$P_{2r}$ is colored~$\{k\}$.
Since this edge coloring is a strong royal $k$-edge coloring of~$C_{2r}$,
it follows that $\sroy(P_{2r})=\sroy(C_{2r})= k$ for  $7 \le r \le 15$, where $k = 4$ if $r = 7$ and $k = 5$ if $8 \le r \le 15$.
Figure~\ref{roy25} shows  the  construction of  a strong royal $4$-edge coloring  of $C_{14}$
from the paths $P_7$ and $P_7^*$.

For each such  path $P_{2r}$ ($7 \le r \le 15$), we construct the path~$P_{2r+1}$ by adding a new vertex~$u_0$ and the edge $u_0u_1$
and coloring the edge $u_0u_1$ by~$\{k\}$, where $k = 4$ if $r = 7$ and $k = 5$ if $8 \le r \le 15$.
Then $u_0$ is colored~$\{k\}$, resulting in  a strong royal $k$-edge coloring of~$P_{2r+1}$
for $7 \le r \le 15$. Next, we repeat this procedure by beginning  with the paths $P_{14}$, $P_{15}$, $\ldots$, $P_{31}$; that is,
we use $P_{14}$  to create a strong royal $5$-edge coloring of~$C_{28}$ (where  $r = 14$) and use
$P_{15}, P_{16},  \ldots, P_{31}$ to create a strong royal $6$-edge coloring of~$C_{2r}$
for $15 \le r \le 31$.  Continued repetition of this procedure gives the desired result for all even cycles.
Therefore, $\sroy(C_{n}) = k$ for all even integers $n \ge 4$ with $2^{k-1} \le n \le 2^k-1$.

\smallskip

\begin{figure}[ht]\centering
\centerline{\input{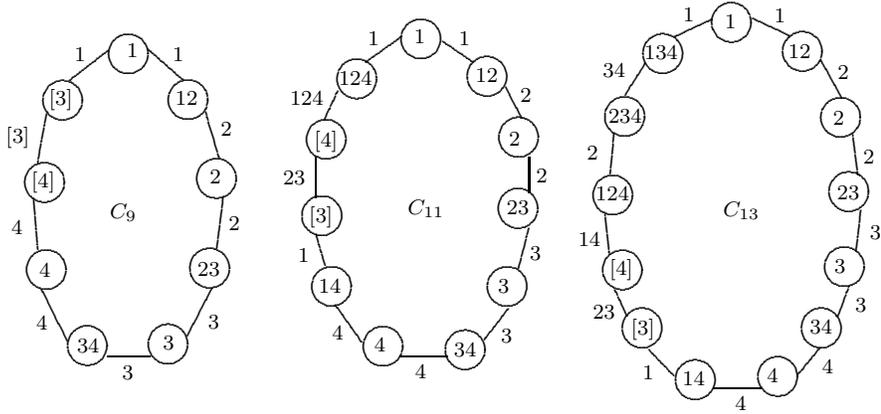}}
\caption{Strong royal $4$-edge colorings of $C_n$  for $n = 9, 11, 13$} \label{roy23}
\end{figure}

{\it Case~$2$.}  {\it $n \ge 9$ is  odd.}   Figure~\ref{roy23}
shows a strong royal 4-edge coloring for each of~$C_9$, $C_{11}$, and~$C_{13}$
and so $\sroy(C_n) = 4$  for $n = 9, 11, 13$.  Thus, we assume that  $n = 2r+1\ge 15$,  where~$r \ge 7$.

\begin{figure}[ht]\centering
\centerline{\input{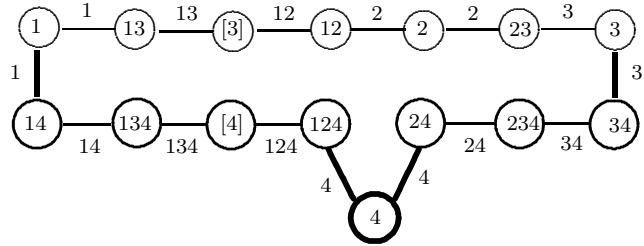}}
\caption{Constructing a strong royal $4$-edge coloring  of $C_{15}$} \label{roy26}
\end{figure}

For each path $P_r$, there is a subpath~$Q=(v_i, v_{i+1}, v_{i+2}, v_{i+3})$,
where $3 \le i <  i+4\le r$ such that
%\bce
$c'(v_{i+1})= \{1, 2\}$,  $c(v_{i+1}v_{i+2}) = \{2\}$, and $c'(v_{i+2})= \{2\}$.
%\ece
From the manner in which  each even cycle $C_{2r}$ was constructed
and a strong royal $k$-edge coloring~$c$ of $C_{2r}$ was defined  in Case~$1$, the path
$Q^*=(u_i, u_{i+1}, u_{i+2}, u_{i+3})$ is a subapth in~$C_{2r}$ such that
%\bce
$c'(u_{i+1})= \{1, 2, k\}$,  $c(u_{i+1}u_{i+2}) = \{2, k\}$, and $c'(u_{i+2})= \{2, k\}$.
%\ece
Furthermore, $c'(x) \ne \{k\}$ for each vertex $x$ of~$C_{2r}$.  We now construct the cycle~$C_{2r+1}$ from
$C_{2r}$ by deleting the edge $u_{i+1}u_{i+2}$ from $C_{2r}$ and adding a new vertex~$u$
along with the two new edges $u_{i+1}u$ and $uu_{i+2}$.
We define an edge coloring~$c$ of $C_{2r+1}$ from the strong royal $k$-edge coloring~$c$ of $C_{2r}$ (as described in Case~$1$)
by assigning the color~$\{k\}$ to the edges $u_{i+1}u$ and $uu_{i+2}$ where the colors of remaining  edges of~$C_{2r+1}$
are the same as in~$C_{2r}$. Thus, $c'(u)=\{k\}$  and $c'(x)$ is the same as in $C_{2r}$ for all other vertices~$x$ of~$C_{2r+1}$.
Figure~\ref{roy26} shows  the  construction of such a strong royal $4$-edge coloring  of~$C_{15}$
from the strong royal $4$-edge coloring  of $C_{14}$  of  Figure~\ref{roy25}.
Since this edge coloring is a strong royal $k$-edge coloring  of $C_{2r+1}$, it follows that
$\sroy(C_{n}) = k$ for all odd integers $n\ge 3$ with $2^{k-1} \le n \le 2^k-1$ with the exception of~$n=3$ and~$n=7$.~\hfill \qed

\medskip
It is therefore a consequence of Theorem~\ref{royalcn}  that~$C_3$ and~$C_7$ are royal-one but all other cycles are royal-zero.

\section{Classes of  Royal-Zero and  Royal-One  Graphs}

In this section we determine some classes of graphs that are royal-zero or royal-one.
For complete graphs, the following result was obtained in~\cite{royc1}.

\bpr                                                 \label{sroyKn}
\ For an integer $n \ge 4$, the complete graph $K_n$ is a royal-zero graph if $n$ is a power of~$2$ and royal-one otherwise.
\epr

We now consider the effect that certain operations can have on graphs that are  royal-zero  or  royal-one.
The  {\it corona} $\corona(G)$ of a graph $G$ is that graph obtained from $G$ by
adding a pendant edge at each vertex of $G$. Thus, if  the order of~$G$ is~$n$, then
the order of~$\corona(G)$ is~$2n$.  The strong royal index of~$\corona(G)$ never exceeds $\sroy(G)$ by more than~1.

\bpr               \label{sroycorona}
\  If $G$ is a connected graph of order~$n \ge 4$, then
  \bce
   $\sroy(\corona(G)) \le \sroy(G) + 1$.
   \ece
Consequently, if $G$ is a royal-zero graph, then so is  $\corona(G)$.
\epr
\pf \quad   Let $V(G)=\{v_1, v_2, \ldots, v_n\}$ and let  $H = \corona(G)$ be obtained from $G$ by adding
the pendant edge $u_iv_i$ at $v_i$ for $1 \le i \le n$.  Suppose that $\sroy(G)= k$. Then there is a
strong royal $k$-edge coloring   $c_G: E(G) \to  \cP^*([k])$ of~$G$.
Define an edge coloring $c_H: E(H) \to  \cP^*([k+1])$ by
$$c_H(e)=\left\{
  \barr{ll}
  c_G(e) \cup \{k+1\} & \mbox{ if $e \in E(G)$}\\[.1cm]
   c'_G(v_i) & \mbox{ if $e= u_iv_i$ for  $1 \le i \le n$.}
   \earr\right.
$$
Then the induced vertex coloring $c'_H$ is given by
\bce
$c_H'(u_i) = c'_G(v_i)$ and  $c'_H(v_i) =  c_G'(v_i) \cup \{k+1\}$ for $1 \le i \le n$.
\ece
Since $c'_H$ is vertex-distinguishing, it follows that  $c_H$ is a strong royal $(k+1)$-edge coloring of~$\corona(G)$ and so
$\sroy(H) \le k+1 =  \sroy(G) + 1$.

If  $G$ is a connected royal-zero graph of order~$n \ge 4$ where $\sroy(G)=k$, say, then  $2^{k-1} \le n \le 2^k-1$.
Since $\corona(G)$  is a connected graph of order~$2n \ge 8$ where $2^{k} \le 2n \le 2^{k+1}-2$,
it follows that $\sroy(\corona(G)) \ge k+1$. On the other hand, there is a  strong royal $(k+1)$-edge coloring of~$\corona(G)$
and so $\sroy(\corona(G))= k+1$, which implies that $\corona(G)$ is royal-zero as well.~\hfill \qed

\medskip

A tree $T$  is called {\it cubic} if every vertex of $T$ that is not an end-vertex has degree~3.
The following result makes use of the proof of Proposition~\ref{sroycorona}.

\bco                                 \label{royalcubiccater}
\   If $T$  is a cubic caterpillar of order at least~$4$, then $T$ is royal-zero.
\eco
\pf \quad Let $T$ be a cubic caterpillar. Since the statement is true if $T$ has four vertices,
we may assume that $T$ has six or more vertices.
For an integer~$n \ge 4$ where $2^{k-1} \le  n  \le 2^k -1$, let
$H =P_n=  (v_1, v_2, \ldots, v_n)$  be a longest path in~$T$, where then $\diam(T) = n-1 \ge 3$ and  the order of~$T$  is~$2n-2$.
As noted earlier, it was shown in~\cite{royc1} that all paths of order~4 or more are royal-zero and so~$\sroy(H)=k$.
Let~$u_iv_i$ be the pendant edges at~$v_i$ for $2 \le i \le n-1$.
We consider two cases, according to whether $2^{k-1} < n \le 2^{k}-1$ or  $n = 2^{k-1}$.
In the first case, we apply the same procedure used in the proof of Proposition~\ref{sroycorona}.
\smallskip

{\it Case~$1$}. {\it $2^{k-1} < n \le 2^{k}-1$.}    Then $2^k \le   2n-2 \le  2^{k+1} -4$.
Thus,  it suffices to show that $\sroy(T) \le  k+1$.
Since $\sroy(H) = k$, there is a strong royal $k$-edge coloring
$c_{H}: E(H) \to \cP^*([k])$. Define an edge coloring $c_T: E(T) \to  \cP^*([k+1])$ by
$$c_T(e)=\left\{
  \barr{ll}
  c_{H}(e) \cup \{k+1\} & \mbox{ if $e \in E(H)$}\\[.1cm]
   c'_{H}(v_i) & \mbox{ if $e= u_iv_i$ for  $2 \le i \le n-1$.}
   \earr\right.
$$
Then the induced vertex coloring $c'_T$ is given by
%\bce
$c_T'(u_i) = c'_{H}(v_i)$ for $2 \le i \le n-1$ and  $c'_T(v_i) =  c_{H}'(v_i) \cup \{k+1\}$ for $1 \le i \le n$.
%\ece
Since $c'_T$ is vertex-distinguishing, it follows that  $c_T$ is a strong royal $(k+1)$-edge coloring of~$T$ and
$\sroy(T) \le k + 1$. Thus,  $T$ is royal-zero.

\smallskip

\begin{figure}[ht]\centering
\centerline{\input{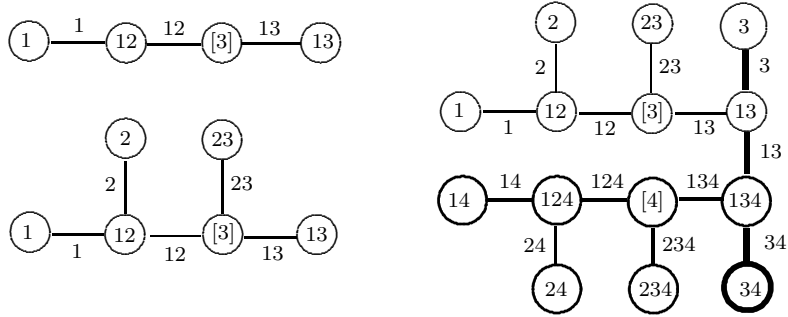}}
\caption{Constructing   strong royal colorings of the cubic caterpillars} \label{roy29}
\end{figure}

{\it Case~$2$}. {\it $n = 2^{k-1}$.} Then $2n-2 = 2^k -2$. Here, we show that $\sroy(T)=\sroy(H)= k$.
First, we consider the case where~$n = 4$ and $k = 3$. A strong royal $3$-edge coloring~$c$ of~$H=P_4=(v_1, v_2, v_3, v_4)$
is shown in Figure~\ref{roy29}, namely $c(v_1v_2)=1$, $c(v_2v_3)=\{1, 2\}$, and $c(v_3v_4)=\{1, 3\}$. Observe that
the induced vertex colors of the vertices of~$H$ are all subsets of~$[3]$ containing~1 and $c'(v_1)=\{1\}$.
The tree~$T$ is constructed from~$H$ by
attaching the pendant edges~$u_2v_2$ and~$u_3v_3$ to~$v_2$ and~$v_3$, respectively. The colors of~$u_iv_i$, $i= 2, 3$,
are  defined by $c(u_iv_i)= c'(v_i)-\{1\}$, which results in a strong royal $3$-edge coloring of~$T$.
In the case where~$n = 8$ and $k = 4$, we begin with the path $H=P_8 =(v_1, v_2, \ldots, v_8)$,  where the edges $v_1v_2$, $v_2v_3$,
$v_3v_4$ of  $P_8$ are colored as in the case when $n = 4$,  and define $c(v_4v_5)=c'(v_4)$ and $c(v_iv_{i+1})=c(v_{8-i}v_{9-i})\cup\{4\}$
for $i = 5, 6, 7$. Here too, each edge color and induced vertex color contains~1 and $c'(v_1)=\{1\}$.
The tree~$T$ in this case is constructed from~$H$ by  attaching the pendant edges~$u_iv_i$ for $2 \le i \le 7$.
The color of $u_iv_i$  is defined by $c(u_iv_i)= c'(v_i)-\{1\}$ for $2 \le i \le 7$, which results in a strong royal $4$-edge coloring of~$T$.
This is illustrated in Figure~\ref{roy29}. Continuing in this manner gives the desired result.~\hfill \qed

\begin{figure}[ht]\centering
\centerline{\input{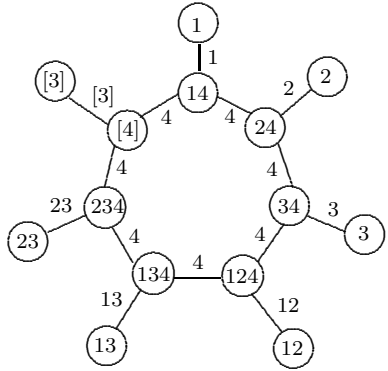}}
\caption{A strong royal 4-edge coloring   of $\corona(C_7)$} \label{roy28}
\end{figure}

As stated in Proposition~\ref{sroycorona},
 if $G$ is a connected graph of order~$4$ or more, then
$\sroy(\corona(G)) \le \sroy(G) + 1$, which implies that
 if $G$ is  a royal-zero graph, then $\corona(G)$ is a royal-zero graph.
On the other hand, it is possible that  $G$ is a royal-one graph and $\corona(G)$ is a royal-zero graph.
By Proposition~\ref{sroyKn},  every complete graph $K_n$ where $n$ is not a power of~$2$ is
a royal-one graph. Thus, if $2^{k-1} +1 \le n \le 2^k-1$ for some integer $k \ge 3$, then
$\sroy(K_n)=k+1$. If one were to assign distinct nonempty subsets of~$[k]$ to the~$n$ pendant edges
of $\corona(K_n)$ and the color~$\{k+1\}$ to the remaining ${n \choose 2}$ edges of~$\corona(K_n)$, then
we have a strong royal $(k+1)$-edge coloring  of $\corona(K_n)$ and so $\sroy(\corona(K_n)) = k+1$.
Therefore,  $\corona(K_n)$ is a royal-zero graph for each integer $n \ge 5$ where $n$ is not a power of~$2$.
For a more interesting example,  Figure~\ref{roy28} shows a  strong royal $4$-edge coloring of~$\corona(C_7)$
 and so  $\sroy(\corona(C_7)) = \sroy(C_7) = 4$ (by Theorem~\ref{royalcn}).
Thus,~$C_7$ is royal-one,  while~$\corona(C_7)$ is royal-zero.

A graph operation somewhat related to the corona of a graph~$G$ is  the {\it  Cartesian product}
of~$G$ with~$K_2$. In fact, we have the following result that corresponds to Proposition~\ref{sroycorona}.

\bpr       \label{sroyGBoxK2}
If $G$ is a connected graph of order~$n \ge 4$, then
  \bce
   $\sroy(G\ \Box \ K_2) \le \sroy(G) + 1$.
   \ece
Consequently, if  $G$ is a royal-zero graph, then $G\ \Box \ K_2$ is a  royal-zero graph.
\epr
\pf \quad   Let $G$  be a  connected graph of order $n \ge 4$ where  $\sroy(G)=k$  for some  positive integer~$k$.
Let  $H=G\ \Box \ K_2$  where $G_1$  and $G_2$ are the two copies of~$G$.
Suppose that  $V(G_1)=\{u_1,u_2,\dots u_n\}$ where $u_i$ is labeled~$v_i$ in~$G_2$.
Thus,  $V(G_2) = \{v_1,v_2,\dots, v_n\}$ and  $E(H)= E(G_1)\cup E(G_2) \cup \{u_iv_i:  1 \le i \le n\}$.
Since  $\sroy(G)=k$, there is  a strong royal $k$-edge coloring $c_{G_1}: E(G_1) \to \cP^*([k])$ of~$G_1$.
Define an edge coloring $c_H: E(H) \to \cP^*([k+1])$ by
 $$c_H(e) =\left\{\barr{ll}
    c_{G_1}(e)  & \mbox{if $e \in E(G_1)$}\\[.1cm]
    c_{G_1}(u_iu_j)\cup \{k+1\}  & \mbox{if $e=v_iv_j \in E(G_2)$ for $1 \le i\ne j \le n$}\\[.1cm]
     c'_{G_1}(u_i) & \mbox{if $ e=u_iv_i$ for $1 \le i \le n$.}
    \earr \right.
    $$
The induced coloring $c'_H: V(H) \to \cP^*([k+1])$ is then given by
$c'_H(u_i)= c'_{G_1}(u_i)$ and $c'_H(v_i) = c'_{G_1}(u_i) \cup\{k+1\}$.
Since $c'_H$ is vertex-distinguishing, it follows that $c'_H$ is a strong royal $(k+1)$-edge coloring of~$H$.
Thus,  $\sroy(H) \le   k+1= \sroy(G)+1$.
Therefore,  if  $G$ is a royal-zero graph, then $G\ \Box \ K_2$ is a  royal-zero graph.~\hfill \qed

\medskip
The {\it hypercube} $Q_k$ is $K_2$ if $k = 1$,
while  for $k \geq 2$, $Q_k$ is defined recursively as the Cartesian product
$Q_{k-1} \ \Box \ K_2$ of $Q_{k-1}$ and $K_2$.
Since $Q_2= C_4$  is royal-zero by Theorem~\ref{royalcn}, the following is a consequence of Proposition~\ref{sroyGBoxK2}.

\bco
\ For each integer $k \ge 2$, the hypercube~$Q_k$ is a royal-zero graph.
\eco

\begin{figure}[ht]\centering
\centerline{\input{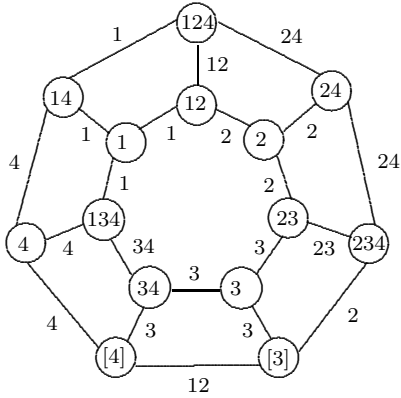}}
\caption{A strong royal 4-edge coloring   of $C_7 \  \Box \ K_2$} \label{roy30}
\end{figure}

As stated in Proposition~\ref{sroyGBoxK2},  if $G$ is a royal-zero graph, then $G\ \Box \ K_2$ is a royal-zero graph.
On the other hand, it is possible that  $G$ is a royal-one graph and $G\ \Box \ K_2$ is a royal-zero graph.
To see an example of this, we return to the $7$-cycle~$C_7$, which we saw (in Theorem~\ref{royalcn})
is a royal-one graph.
Figure~\ref{roy30} shows a  strong royal $4$-edge coloring of~$C_7 \  \Box \ K_2$
and so  $\sroy(C_7) = \sroy(C_7 \  \Box \ K_2) = 4$ .
Thus,~$C_7$ is royal-one,  while~$C_7 \ \Box \ K_2$ is royal-zero. As mentioned in Proposition~\ref{sroyKn},
the complete graphs~$K_5$ and~$K_6$ are royal-one graphs. For these two
graphs~$G$, the graphs~$G\ \Box \ K_2$ are royal-zero; that is,
$\sroy(K_5 \ \Box \ K_2) = \sroy(K_6 \ \Box \ K_2) = 4$.  A   strong royal $4$-edge coloring~$c$ of~$H=K_6 \  \Box \ K_2$
can be defined as follows. Let $H_1$ and $H_2$ be two copies of~$K_6$ in $H$,
where $V(H_1)=\{u_1,u_2,\dots u_6\}$ and $V(H_2)= \{v_1,v_2,\dots v_6\}$ such that $u_iv_i \in E(H)$.
First, we define the vertex-distinguishing  coloring $c':  V(H) \to \cP^*([4])$ by
\bce
$c'(u_1)=\{1, 4\}$, $c'(u_2)=\{1\}$, $c'(u_3)=\{1, 2, 4\}$,

$c'(u_4)=\{1, 2, 3\}$, $c'(u_5)=\{1, 3\}$, $c'(u_6)=\{1, 2\}$,

$c'(v_1)=\{4\}$, $c'(v_2)=\{1, 3, 4\}$, $c'(v_3)=[4]$,

$c'(v_4)=\{2, 4\}$, $c'(v_5)=\{3, 4\}$, $c'(v_6)=\{2, 3, 4\}$.
\ece
The edge coloring $c:  V(H) \to \cP^*([4])$ is then defined  by $c(xy) = c'(x) \cap c'(y)$  for each edge $xy \in E(H)$.
Since~$c'$ is  the induced vertex coloring of~$c$, it follows that $c$ is a  strong royal $4$-edge coloring  of~$H$.
Thus, $H$ is  royal-zero.

As noted in Proposition~\ref{sroyKn},
the complete graph $K_7$ is also  a  royal-one graph.  However, $H=K_7\ \Box \ K_2$  is   royal-one as well.
That there is a  strong royal $5$-edge coloring  of~$H$ is  straightforward. To show that
$\sroy(K_7 \ \Box \ K_2) = 5$, however, it is necessary to show that there is no strong royal $4$-edge coloring  of~$H$,
for assume that such an edge coloring $c$ of~$H$ exists.  Since the order of~$H$ is 14, the induced
vertex colors of~$H$ must consist of 14 elements of~$\cP^*([4])$. In particular, at least three of the four singleton subsets
of~$[4]$ must be vertex colors of~$H$. Suppose that~$H_1$ and~$H_2$ are the two copies of~$K_7$ in the construction of~$H$.
Therefore, at least one of~$H_1$ and~$H_2$ has at least two singleton subsets as its vertex colors, say $c'(u_1)=\{1\}$ and
$c'(u_2)=\{2\}$ where $u_1, u_2 \in V(H_1)$, which is impossible since~$u_1$ and~$u_2$ are adjacent. Hence, $\sroy(K_7 \ \Box \ K_2) = 5$.

\section{Conditions  for  Royal-One Graphs}         \label{se2}

\begin{figure}[ht]\centering
\centerline{\input{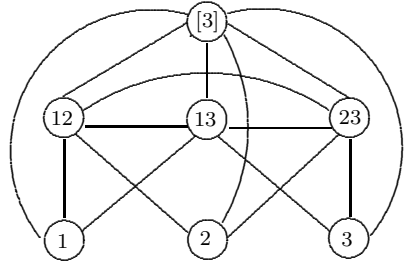}}
\caption{The graph $G_3$ of order $7=2^3 -1$ and size~$m_3=15$} \label{roy18}
\end{figure}

We have seen that  many graphs are royal-zero graphs.
We now  present a sufficient condition for a connected graph~$G$ of order~$n \ge 4$ to be a royal-one graph.
Let $k$ be  the unique integer  such that $2^{k-1} \le n \le 2^k-1$.
A graph~$G_k$ of order~$2^k-1$  is now constructed as follows.
The vertices of~$G_k$  are labeled with the
$2^k-1$ distinct  elements of  $\cP^*([k])$.
For each  vertex~$v$ of~$G_k$, let $\ell(v)$ denote  its  label. Thus,
   %\bce
   $\{\ell(v):  \ v \in V(G_k)\} = \cP^*([k]).$
  % \ece
Two vertices $u$ and $v$ of~$G_k$ are adjacent in $G_k$  if and only if $\ell(u) \cap \ell(v) \ne \emptyset$.
The vertex set $V(G_k)$ is partitioned into $k$ subsets $V_1, V_2, \ldots, V_k$
where $V_i = \{v\in V(G_k) : |\ell(v)| = i\}$ for $1 \le i \le k$. Therefore, $G_k[V_k] =K_1$ and $G_k[V_1] =\ol{K}_k$ is empty.
If $k =2p+ 1$ is odd, then $G_k\left[V_{p+1} \cup V_{p+2} \cup  \cdots  \cup V_k\right] = K_{2^{k-1}}$.
If $k = 2p$ is even, then let $V_{p}'$ be the subset  consisting of  those elements~$S$ in $V_{p}$
for which $1 \in S$. Then $|V_{p}'| = \frac{1}{2}{k \choose p}$ and
$G_k\left[V'_{p} \cup V_{p+1} \cup V_{p+2} \cup  \cdots  \cup V_k\right] = K_{2^{k-1}}$.
Let  $m_k$  be the size of~$G_k$.
The graph $G_3$ of order $7=2^3 -1$ has size~$m_3= 15$ and is shown in Figure~\ref{roy18}.

There is an immediate  condition under which a connected graph cannot be a royal-zero graph.
As we mentioned earlier, it  was shown in~\cite{BDK2017} that    $G$ is a connected graph of order~$n \ge 4$  where $2^{k-1} \le n \le 2^k-1$, then $\sroy(G) \le k+2$.

\bobs                     \label{Gkmkde}
\ Let  $G$ be  a  connected graph of order $n \ge 4$ and size~$m$ where $2^{k-1} \le n \le 2^k-1$ for an integer~$k$.
If $G$ is not a subgraph of the graph~$G_k$, then  either $\sroy(G) = k+1$ or $\sroy(G) = k+2$,  and so $G$ is not a royal-zero graph.
Consequently, if  $m \ge m_k+1$, then $G$ is  not a royal-zero graph.
\eobs

\begin{figure}[ht]\centering
\centerline{\input{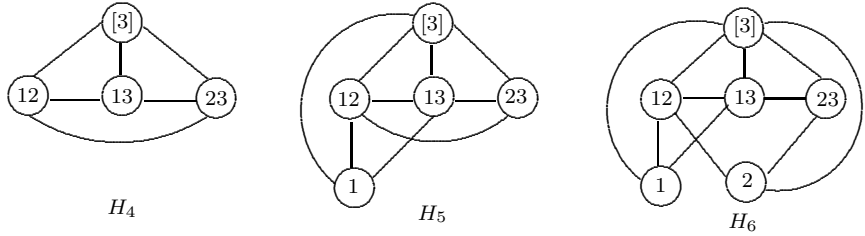}}
\caption{Subgraphs of $G_3$} \label{roy20}
\end{figure}

Since  $\sroy(T)= 3$ for each  tree~$T$ of order~$n$ where $4 \le n\le 7$,  it follows by
Observation~\ref{royalspanH} that   if $G$ is a connected graph of order~$n$ where $4 \le n\le 7$,
then  $\sroy(G)$ is either~3 or~4. If $G$ is a connected graph of order~7 that is not isomorphic to a subgraph of~$G_3$ of Figure~\ref{roy18},
then $\sroy(G) \ne 3$ and so $\sroy(G)=4$. Since the size of $G_3$ is~15, it follows that
if $G$ is a connected graph of order~7 with size at least~$16$,   then  $\sroy(G)=4$.
Figure~\ref{roy20} shows the graphs $H_4, H_5$, and $H_6$ of order~4, 5, and~6, respectively, of greatest size that are subgraphs of~$G_3$.
For each  graph~$H_i$ where $i = 4, 5, 6$,  if  every edge $uv$ of~$H_i$ is  assigned the color~$c(uv) = \ell(u) \cap \ell(v)$, then
$c'(v)=\bigcup_{e \in E_{H_i}(v)} c(e)  =\ell(v)$, resulting in a strong royal $3$-edge coloring of~$H_i$.
Hence, $\sroy(H_i) = 3$ for $i = 4, 5, 6$.
The graph~$H_4=K_4$, while $H_5$ has size~9 and~$H_6$ has size~12.
So, if $G$ is a connected graph of order~5 whose  size  is at least~10 (that is, $G = K_5$) or
 if $G$ is a connected graph of order~6  whose  size  is at least~13, then  $\sroy(G)=4$.

\begin{figure}[ht]\centering
\centerline{\input{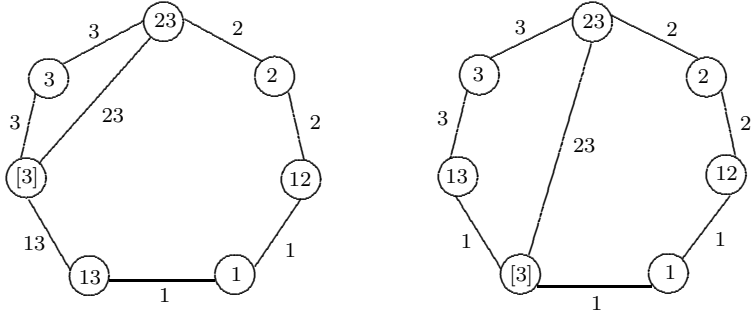}}
\caption{Showing   that $\sroy(C_7 + e)=3$ for each $e \notin E(C_7)$} \label{roy22}
\end{figure}

By Observation~\ref{Gkmkde},    if $G$ is  a  connected graph of order $n \ge 4$ and size~$m$ where $2^{k-1} \le n \le 2^k-1$
such that $m  >  m_k$, which implies that  $G \not\sbe G_k$, then $\sroy(G) \ge k+1$.
In fact, if $G$ possesses any property that implies that $G \not\sbe G_k$, then $\sroy(G) \ge k+1$.
For example, if the order of~$G$ is~$n = 2^k-1$ and $\de(G) \ge \de(G_k)+1$ or $G$ has more than one vertex of degree~$n-1$, then
$\sroy(G) \ge k+1$.
On the other hand,  even though $C_7 \sbe G_3$ (where $n = 2^3-1$ and $k = 3$),
$|E(C_7)|= 7 < m_3$,  and  $\de(C_7) < \de(G_3)$, we saw that $\sroy(C_7) = 4 = k+1$.
Furthermore, for every chord  $e$ of~$C_7$, $\sroy(C_7+e)=3$ (see Figure~\ref{roy22}).
Consequently, even though one might suspect that
$\sroy(G + uv) \ge \sroy(G)$ for every connected graph~$G$ and every pair $u, v$ of nonadjacent vertices of~$G$, such is not  the case.

What we have seen is that if $G$ is a connected graph  of order $n \ge 4$ where $2^{k-1} \le n \le 2^k-1$ having a sufficiently large size,
then $\sroy(G) \ne k$. However, if $G$  is a connected graph  of order $n \ge 4$ where $2^{k-1} \le n \le 2^k-1$ having a small size, then we are not
guaranteed that $\sroy(G) =  k$.  Indeed, even the strong royal index of trees is in doubt.

If Conjecture~\ref{royalgraph2.1} is true,  then for every connected graph~$G$ of order $n \ge 4$ where $2^{k-1} \le n \le 2^k-1$,
either $\sroy(G)=k$ or $\sroy(G)=k+1$.  In order to present a sufficient condition for
$\sroy(G)\ne k$ in terms of the size and minimum degree of~$G$, we describe an expression for the
size~$m_k$ of the graph~$G_k$ (as it is easier in general to compare two numbers than to determine whether
a graph contains a subgraph isomorphic to a given graph).

Recall that we label the $2^k-1$  vertices of~$G_k$  with the distinct  elements of  $\cP^*([k])$.
The label of each vertex $v$ of~$G_k$ is denoted by~$\ell(v)$ and so $\{\ell(v):  \ v \in V(G_k)\} = \cP^*([k]).$
Let  $\{V_1, V_2, \ldots, V_k\}$ be the partition of of~$V(G_k)$ described earlier, where  then
$V_i = \{v\in V(G_k) : |\ell(v)| = i\}$ for $1 \le i \le k$.   Let $v \in V_i$ for some integer~$i$ with $1 \le i \le k$.
Then $\ell(v) = S$ is some $i$-element subset of~$[k]$. There are~$2^i-1$ nonempty subsets of~$S$ and~$2^{k-i}$
subsets of~$[k]-S$. For each nonempty subset~$S'$ of~$S$ and each subset~$T$ of~$[k]-S$, the vertex~$v$ is adjacent to
that vertex~$w$ of~$G_k$ for which~$\ell(w)=S'\cup T$. Since~$v$ is not adjacent to itself, however,
it follows that $\deg_{G_k} v = (2^i-1)2^{k-i} - 1$.  Furthermore,  there are ${k \choose i}$ vertices in~$V_i$ for $1 \le i \le k$.
Therefore,
\bea
m_k &=& \frac{1}{2}\sum_{i=1}^k {k \choose i}  \left[(2^i-1)2^{k-i} - 1\right] = \frac{1}{2}\sum_{i=1}^k {k \choose i}  (2^k - 2^{k-i} -1)\\[.1cm]
        &=& \frac{1}{2}\left[\sum_{i=1}^k {k \choose i}2^k -     \sum_{i=1}^k {k \choose i}2^{k-i} -     \sum_{i=1}^k {k \choose i}\right]\\[.1cm]
        &=& \frac{1}{2}\left[2^k \sum_{i=1}^k {k \choose i} - 2^k  \sum_{i=1}^k {k \choose i}\left(\frac{1}{2}\right)^i -     \sum_{i=1}^k {k \choose i}\right]\\[.1cm]
       &= & \frac{1}{2}\left\{ 2^k(2^k-1) - 2^k\left[\left(1+\frac{1}{2}\right)^k -1\right] - (2^k-1)\right\} \\[.1cm]
        &=&   \frac{1}{2}(4^k-3^k-2^k+1).
        \eea
In particular,  if $k = 3$, then the size of~$G_3$ is $m_3= 15$, as we saw in Figure~\ref{roy18}.

\bpr
\ Let $G$ be a graph of order $n \ge 4$ and size~$m$ where $2^{k-1} \le n \le 2^k-1$ for some integer~$k\ge 3$.
If $m >  \frac{1}{2}(4^k-3^k-2^k+1)$, then   either $\sroy(G) = k+1$ or $\sroy(G) = k+2$,  and so $G$ is not a royal-zero graph.
\epr

For each integer $k \ge 3$, the minimum degree~$\de(G_k)$ of the graph~$G_k$ is~$2^{k-1}-1$. Consequently,
if  $G$ is a graph of order $n \ge 4$ and size~$m$ where $2^{k-1} \le n \le 2^k-1$
for which $\de(G) \ge 2^{k-1}$, then it may occur that $m < m_k$ but yet $G$ is not a subgraph of~$G_k$,  and so
(by  Observation~\ref{Gkmkde}) $\sroy(G)\ge  k + 1$.  However, in this case, more can be said.
It is useful to recall that every path~$P_n$ for $n \ge 4$ is royal-zero~(see~\cite{BDK2017, royc1}).

\bpr       \label{sroydeG}
\ Let  $G$  be a connected graph  of order $n\ge 4$ where  $2^{k-1} \le n \le 2^k-1$ for some integer $k \ge 2$.
If $\de(G)\ge 2^{k-1}$, then  $\sroy(G) =   k + 1$ and so $G$ is a royal-one graph.
\epr
\pf \quad   We have already observed that $\sroy(G)\ge  k + 1$ for such a graph.
Since $\de(G)\ge 2^{k-1}$ and $n \le 2^k-1$, it follows that
$\de(G) \ge (n+1)/2$ and therefore $G$ has a Hamiltonian path (in fact, a Hamiltonian cycle).  Since
$\sroy(P_n)= k$ for every path~$P_n$ of order~$n$,
it follows by Observation~\ref{royalspanH}  that
$\sroy(G) \le k+1$ and so $\sroy(G)= k+1$.~\hfill \qed

\medskip

We have seen that both~$K_7$ and~$C_7$ (a spanning subgraph,  or factor,  of~$K_7$)
are royal-one graphs. The complement~$\ol{C}_7$ of~$C_7$ is a $4$-regular graph of order~7 and so it is not a subgraph
of the graph~$G_3$ shown in Figure~\ref{roy18}.  Hence,    $\ol{C}_7$ is also  a royal-one graph. The size of~$\ol{C}_7$
is 14 which is less than the size~15 of~$G_3$ (the graph of order~7 having the maximum size that is royal-zero).
This brings up the problem of determining for each integer~$n \ge 3$, the  minimum  size of a graph of order~$n$ that is royal-one.
Of course, the minimum size is~7 when $n = 7$.

The graph~$\ol{C}_7$ can itself be factored into two copies of~$C_7$. Therefore, the royal-one graph~$K_7$ can be factored into three
royal-one graphs. However, $K_7$ can also  be factored into three graphs satisfying any of the following:
(1) all three factors are royal-zero,
(2) exactly two factors are royal-zero,
(3)    exactly one factor is royal-zero.
Consequently, there is a host of additional problems that arise with strong royal colorings of graphs.

\end{document}